\documentclass[a4paper,reqno]{amsart}
\usepackage{amssymb, amsmath, amscd}
\usepackage{color}

\def\gronk{\noalign{\hrule}\vphantom{\vrule height 13pt\frac1{4_{A_A}}}}
\newtheorem{theorem}{Theorem}[section]
\newtheorem{definition}[theorem]{Definition}

\newtheorem{example}[theorem]{Example}
\newtheorem{remark}[theorem]{Remark}

\def\WW{\mathcal{W}}\def\WWW{\tilde{\mathcal{W}}}

\def\GL{\operatorname{GL}}
\def\mapright#1{\smash{\mathop{\longrightarrow}\limits\sp{#1}}}
\makeatletter
 \@addtoreset{equation}{section}
\makeatother
\begin{document}
\title[Geometric realizations]
 {Geometric realizations of curvature}
\author{M.  Brozos-V\'azquez, P. Gilkey, and S. Nik\v cevi\'c}
\address{MB: Department of Mathematics, University of A Coru\~na, Spain}
\email{mbrozos@udc.es}
\begin{address}{PG: Mathematics Department, University of Oregon,
Eugene Or 97403 USA.}
\end{address}
\email{gilkey@uoregon.edu}
\begin{address}
{SN: Mathematical Institute, Sanu,
Knez Mihailova 35, p.p. 367,
11001 Belgrade,
Serbia}\end{address}
\email{stanan@mi.sanu.ac.rs}

\begin{abstract}{We study geometric realization questions of curvature
in the affine, Riemannian, almost Hermitian, almost para Hermitian, almost hyper Hermitian, almost hyper para Hermitian,
Hermitian, and para Hermitian settings. We also express questions in Ivanov--Petrova geometry, Osserman geometry, and
curvature homogeneity in terms of geometric realizations.}\end{abstract}
\maketitle
\centerline{This paper is dedicated to Professor Sekigawa}
\section{Introduction}
A central area of study in Differential Geometry is the examination of the relationship between
purely algebraic properties of the Riemann curvature tensor and the underlying geometric properties of the manifold.
Many authors have worked in this area in recent years. Nevertheless, many fundamental questions remain
unanswered.

It is frequently convenient to work first purely algebraically and pass later to the geometric setting; many questions in differential
geometry can be phrased as problems involving the geometric realization of curvature. Here is a brief outline to this paper. In Section
\ref{sect-2-Affine}, we study the affine setting, in Section \ref{sect-3-RiemannianGeometry}, we study pseudo Riemannian
geometry, and in Section
\ref{S-4-AffineRiemann}, we combine these two structures and consider realization problems related to affine geometry where
the additional structure of a pseudo Riemannian metric is present. In Section \ref{S5-almostHerm}, we turn our attention to
almost Hermitian and almost para Hermitian geometry and study the scalar curvature and the $\star$-scalar curvature. In
Section
\ref{S6-hyperHerm}, we examine similar questions in hyper almost Hermitian geometry and hyper almost para Hermitian
geometry. In Section \ref{S7-HermGeo}, we study realization questions which arise when the structures in question are to be
integrable.  In Sections \ref{sect-9-IvanovPetrovaGeometry} and \ref{sect-10-OssermanGeometry}, we
discuss Ivanov--Petrova geometry and Osserman geometry, respectively. In Section \ref{sect-11-Curvaturehomogeneity}, we
present questions of curvature homogeneity.

The decomposition of the appropriate space of tensors into
irreducible modules under the appropriate structure group is central to our investigation and we review the appropriate results in each
section. The results in Sections
\ref{sect-2-Affine}-\ref{S7-HermGeo}, although they involve non-linear analysis, are closely tied to the representation theory of the
appropriate group. In contrast, the results of Sections \ref{sect-9-IvanovPetrovaGeometry}-\ref{sect-11-Curvaturehomogeneity} are
non-linear in their very formulation since one is studying orbit spaces under the structure group which are not linear subspaces.
Throughout this paper, we shall let $M$ be a smooth manifold of
dimension $m\ge4$; there are similar results in the $2$ dimensional and $3$ dimensional settings. We shall let $\nabla$ be a
torsion free connection on the tangent bundle of $M$. Let $g$ be a pseudo Riemannian metric of signature $(p,q)$ on $M$ and
let
$\mathcal{M}:=(M,g)$ be the associated pseudo Riemannian manifold.

\section{Affine Geometry}\label{sect-2-Affine}
We refer to \cite{GNS09,SSV91} for further information concerning affine geometry. An {\it affine manifold} is a pair $(M,\nabla)$
where $M$ is a smooth manifold and where $\nabla$ is a torsion free connection on $M$. The associated curvature operator
$\mathcal{R}$ is defined by setting:
$$\mathcal{R}(x,y):=\nabla_x\nabla_y-\nabla_y\nabla_x-\nabla_{[x,y]}\,.$$
This $(1,3)$ tensor satisfies the identities:
\begin{equation}\label{eqn-2.a}
\mathcal{R}(x,y)=-\mathcal{R}(y,x),\quad
\mathcal{R}(x,y)z+\mathcal{R}(y,z)x+\mathcal{R}(z,x)y=0\,.
\end{equation}

Let $V$ be a vector space of dimension $m$. A tensor
$\mathcal{A}\in \otimes^2V^*\otimes\text{End}(V)$ satisfying the symmetries given in Equation (\ref{eqn-2.a})
is called an {\it affine algebraic curvature operator}; let $\mathfrak{A}(V)\subset \otimes^2V^*\otimes\text{End}(V)$ be the
subspace of all such operators. An affine curvature operator $\mathcal{A}\in\mathfrak{A}(V)$ is said to be {\it
geometrically realizable} if there exists an affine manifold
$(M,\nabla)$, if there exists a point $P$ of $M$ (which is called the realizing point), and if there
exists an isomorphism
$\phi:V\rightarrow T_PM$ so that
$\phi^*\mathcal{R}_P=\mathcal{A}$.  In either the algebraic or the geometric setting, one defines the {\it Ricci
tensor}
$\rho$ by:
$$\rho(x,y):=\operatorname{Tr}\{z\rightarrow\mathcal{A}(z,x)y\}\,.$$
\subsection{The decomposition of $\mathfrak{A}(V)$ as a $\operatorname{GL}(V)$ module}
The action of the general linear group $\operatorname{GL}(V)$ on the vector space of affine algebraic curvature operators
$\mathfrak{A}(V)$ by pullback is not irreducible, but decomposes as the direct sum of irreducible modules.
The decomposition $V^*\otimes V^*=\Lambda^2(V^*)\oplus S^2(V^*)$ is a $\operatorname{GL}(V)$ equivariant decomposition of $V^*\otimes V^*$
into irreducible $\operatorname{GL}(V)$ modules; we let $\rho_a$ and $\rho_s$ be the components in $\Lambda^2(V^*)$ and $S^2(V^*)$,
respectively, where
$$
\rho_a(x,y):=\textstyle\frac12\{\rho(x,y)-\rho(y,x)\}\quad\text{and}\quad
\rho_s(x,y):=\textstyle\frac12\{\rho(x,y)+\rho(y,x)\}\,.
$$
One has the following result of Bokan \cite{Bo90} and of Strichartz \cite{S88}:

\begin{theorem}\label{thm-2.1}
Let $\dim(V)\ge4$. The Ricci tensor $\rho$ yields a $\operatorname{GL}(V)$ equivariant short exact sequence
$$0\rightarrow\ker(\rho)\rightarrow\mathfrak{A}(V)\mapright{\rho_a\oplus\rho_s}
\Lambda^2(V^*)\oplus S^2(V^*)\rightarrow0$$
which is equivariantly split by the map $\sigma$ where
\begin{eqnarray*}
&&\{\sigma\rho_a\}(x,y)z=\textstyle\frac{-1}{1+m}\{2\rho_a(x,y)z+\rho_a(x,z)y-\rho_a(y,z)x\},\\
&&\{\sigma\rho_s\}(x,y)z=\textstyle\frac1{1-m}\{\rho_s(x,z)y-\rho_s(y,z)x\}\,.
\end{eqnarray*}
One has a direct sum
decomposition of $\mathfrak{A}(V)$ into irreducible $\operatorname{GL}(V)$ modules:
$$\mathfrak{A}(V)=\ker(\rho)\oplus\Lambda^2(V^*)\oplus S^2(V^*)\,.$$\end{theorem}

We note for the sake of completeness that:
$$\begin{array}{|l|l|}\gronk
\dim\{\ker(\rho)\cap\mathfrak{A}(V)\}=\frac{m^2(m^2-4)}3&
  \dim\{\Lambda^2(V^*)\}=\frac{m(m-1)}2\\\gronk
\dim\{S^2(V^*)\}=\textstyle\frac{m(m+1)}2&\dim\{\mathfrak{A}(V)\}=\textstyle\frac{m^2(m^2-1)}3
\\\noalign{\hrule}
\end{array}$$

\begin{definition}\label{defn-2.2}
\rm Let $\mathcal{A}\in\mathfrak{A}(V)$.\begin{enumerate}
\item $\mathcal{A}$ is {\it Ricci symmetric} if and only if  $\rho\in S^2(V^*)$, i.e. $\rho_a=0$.
\item $\mathcal{A}$ is {\it Ricci anti-symmetric} if and only if $\rho\in\Lambda^2(V^*)$, i.e. $\rho_s=0$.
\item $\mathcal{A}$ is {\it Ricci flat} if and only if $\rho=0$.
\item The {\it Weyl projective curvature operator} $\mathcal{P}$ is the component of $\mathcal{A}$ in
$\ker(\rho)$, i.e. $\mathcal{P}:=\mathcal{A}-\sigma\rho\in\ker\rho$. $\mathcal{A}$ is {\it projectively flat} if and only if
$\mathcal{P}=0$.
\item $\mathcal{A}$ is {\it flat} if and only if $\mathcal{A}=0$, i.e. $\mathcal{A}$ is both projectively flat
and Ricci flat.
\end{enumerate}\end{definition}

\subsection{Equiaffine geometry} Ricci symmetric torsion free connections are often called {\it equiaffine}; they play a central role in
many settings -- see, for example, the discussion in \cite{BGNS06, BDS03, M03, MS99, PSS94}.
The following result is well known \cite{SS-62} and motivates their investigation:

\begin{theorem}\label{thm-2.3}
Let $(M,\nabla)$ be an affine manifold. The following assertions are equivalent:
\begin{enumerate}
\item $\operatorname{Tr}(\mathcal{R})=0$.
\item $\nabla$ is Ricci symmetric.
\item $\nabla$ locally admits a parallel volume form.
\end{enumerate}\end{theorem}

\subsection{Geometric realizability of affine algebraic curvature operators I} Theorem \ref{thm-2.1} gives rise to
additional geometric realizability questions; the decomposition of
$\mathfrak{A}(V)$ as a
$\GL(V)$ module has 3 components so there are 8 natural geometric realization questions which are $\GL(V)$
equivariant. We refer to the discussion in
\cite{GN08,GNW09} for the proof of the following result which shows, in particular, that the symmetries of Equation (\ref{eqn-2.a})
generate the universal symmetries of the curvature operator of a torsion free connection:

\begin{theorem}\label{thm-2.4}
\ \begin{enumerate}
\item Any affine algebraic curvature operator can be geometrically realized by an affine manifold.
\item  Any Ricci symmetric affine algebraic curvature operator can be geometrically realized by a Ricci symmetric
affine manifold.
\item  Any Ricci anti-symmetric affine algebraic curvature operator can be geometrically realized by a Ricci anti-symmetric
affine manifold.
\item  Any Ricci flat affine algebraic curvature operator can be geometrically realized by a Ricci flat
affine manifold.
\item Any projectively flat affine algebraic curvature operator can be geometrically realized by a projectively flat affine manifold.
\item Any projectively flat Ricci symmetric affine algebraic curvature operator can be geometrically realized by a projectively flat Ricci
symmetric affine manifold.
\item A projectively flat Ricci anti-symmetric affine algebraic curvature operator which is not flat can not be geometrically
realized by a projectively flat, Ricci anti-symmetric affine manifold.
\item If $\mathcal{A}$ is flat, then $\mathcal{A}$ is geometrically realized by a flat affine manifold.
\end{enumerate}
\end{theorem}

These geometric realizability results can be summarized in the following table; the non-zero components of
$\mathcal{A}$ are indicated by
$\star$.
$$\begin{array}{|c|c|c|r|c|c|c|r|}\noalign{\hrule}
\ker(\rho)&S^2(V^*)&\Lambda^2(V^*)&&\ker(\rho)&S^2(V^*)&\Lambda^2(V^*)&\\
\noalign{\hrule}\star&\star&\star&\text{yes}&0&\star&\star&\text{yes}\\
\noalign{\hrule}\star&\star&0&\text{yes}&0&\star&0&\text{yes}\\
\noalign{\hrule}\star&0&\star&\text{yes}&0&0&\star&\text{no}\\
\noalign{\hrule}\star&0&0&\text{yes}&0&0&0&\text{yes}\\\noalign{\hrule}
\end{array}$$

\section{Pseudo Riemannian Geometry}\label{sect-3-RiemannianGeometry}
Let $V$ be a finite dimensional real vector space of dimension $m$. One says that $A\in\otimes^4(V^*)$ is an {\it algebraic
curvature tensor} on $V$ if $A$ satisfies the symmetries of the Riemann curvature tensor:
\begin{equation}\label{eqn-3.a}
\begin{array}{l}
A(x,y,z,w)=-A(y,x,z,w)=A(z,w,x,y),\\
A(x,y,z,w)+A(y,z,x,w)+A(z,x,y,w)=0\,.\end{array}
\end{equation}
Let $\mathfrak{R}(V)$ be the space of all such $4$-tensors; note that $\langle\cdot,\cdot\rangle$ induces a non-degenerate
innerproduct on $\mathfrak{R}(V)$. We say that
$\mathfrak{M}:=(V,\langle\cdot,\cdot\rangle,A)$ is a {\it curvature model} if
$A\in\mathfrak{R}(V)$ and if
$\langle\cdot,\cdot\rangle$ is a non-degenerate symmetric bilinear form of signature $(p,q)$ on $V$. $\mathfrak{M}$ is said to be {\it
Riemannian} if $p=0$ and {\it Lorentzian} if $p=1$. Two curvature models
${\mathfrak{M}_1=(V_1,\langle\cdot,\cdot\rangle_1,A_1})$ and
${\mathfrak{M}_2=(V_2,\langle\cdot,\cdot\rangle_2,A_2})$ are said to be isomorphic, and one writes $\mathfrak{M}_1\approx\mathfrak{M}_2$, if there
is an isomorphism
$\phi:V_1\rightarrow V_2$ so that
$$
\phi^*\langle\cdot,\cdot\rangle_2=\langle\cdot,\cdot\rangle_1\quad\text{and}\quad
\phi^*A_2=A_1\,.
$$

\subsection{The decomposition of $\mathfrak{R}(V)$ as an $O(V,\langle\cdot,\cdot\rangle)$ module} If $\mathfrak{M}$ is a
curvature model, then the associated orthogonal group $O(V,\langle\cdot,\cdot\rangle)$ is defined by setting:
$$O(V,\langle\cdot,\cdot\rangle):=\{T\in\operatorname{GL}(V):T^*\langle\cdot,\cdot\rangle=\langle\cdot,\cdot\rangle\}\,.$$
Let
$\varepsilon_{ij}$ and
$A_{ijkl}$ be the components of
$\langle\cdot,\cdot\rangle$ and of
$A$ relative to a basis $\{e_i\}$ for $V$:
$$\varepsilon_{ij}:=\langle e_i,e_j\rangle\quad\text{and}\quad A_{ijkl}:=A(e_i,e_j,e_k,e_l)\,.$$
Let $\varepsilon^{ij}$ be the
inverse matrix. Adopt the {\it Einstein convention} and sum over repeated indices. The components of the {\it Ricci tensor}
$\rho$ and the {\it scalar curvature} $\tau$ are:
$$\rho_{il}:=\varepsilon^{jk}A_{ijkl}\qquad\text{and}\qquad\tau:=\varepsilon^{il}\rho_{il}\,.$$

Let $S_0^2(V^*,\langle\cdot,\cdot\rangle)\subset S^2(V^*)$ be the subspace of trace free symmetric $2$-tensors and let
$\rho_0:=\rho-\frac1m\tau\langle\cdot,\cdot\rangle$ be the trace free Ricci tensor. We refer to Singer and Thorpe \cite{ST} for:
\begin{theorem} Let $\dim(V)\ge4$. There is an $O(V,\langle\cdot,\cdot\rangle)$ equivariant short exact sequence
$$0\rightarrow\ker(\rho)\rightarrow\mathfrak{R}(V)\mapright{\rho_0\oplus\tau}
S_0^2(V^*,\langle\cdot,\cdot\rangle)\oplus\mathbb{R}\rightarrow0$$ which is equivariantly split by the map $\sigma$ where
\begin{eqnarray*}
\sigma(\rho)(x,y,z,w):&=&{\textstyle\frac1{m-2}}\{\rho(x,w)\langle y,z\rangle+\langle x,w\rangle\rho(y,z)\}\\
&-&{\textstyle\frac1{m-2}}\{\rho(x,z)\langle y,w\rangle+\langle x,z\rangle\rho(y,w)\}\\
&-&{\textstyle\frac{\tau}{(m-1)(m-2)}}\{\langle x,w\rangle\langle y,z\rangle-\langle x,z\rangle\langle y,w\rangle\}\,.
\end{eqnarray*}
One has an orthogonal decomposition of $\mathfrak{R}(V)$ into irreducible
$O(V,\langle\cdot,\cdot\rangle)$ modules:
$$\mathfrak{R}(V)=\ker(\rho)\oplus S_0^2(V^*,\langle\cdot,\cdot\rangle)\oplus\mathbb{R}\,.$$
\end{theorem}

We note for the sake of completeness that:
$$\begin{array}{|l|l|}\gronk
\dim\{\mathbb{R}\}=1&
\dim\{\ker(\rho)\cap\mathfrak{R}(V)\}=\textstyle\frac{m(m+1)(m+2)(m-3)}{12}\\\gronk
\dim\{\mathfrak{R}(V)\}=\textstyle\frac{m^2(m^2-1)}{12}&
\dim\{S_0^2(V^*,\langle\cdot,\cdot\rangle)\}=\frac{(m-1)(m+2)}2
\\\noalign{\hrule}
\end{array}$$

\subsection{Geometric realizability of algebraic curvature tensors} Assume given
a pseudo Riemannian manifold $\mathcal{M}:=(M,g)$ of signature
$(p,q)$. Let $\nabla$ be the Levi-Civita connection of $\mathcal{M}$ and let
$R\in\otimes^4T^*M$ be the curvature tensor:
$$R(x,y,z,w)=g(\mathcal{R}(x,y)z,w)\,.$$
Let $\mathfrak{M}=(V,\langle\cdot,\cdot\rangle,A)$ be a curvature model. We say that $\mathfrak{M}$ is geometrically realizable if thereexists a pseudo Riemannian manifold $\mathcal{M}$, if there exists a point $P$ of $M$, and if there exists an isomorphism
$\phi:V\rightarrow T_PM$ so that
$$\phi^*g_P=\langle\cdot,\cdot\rangle\quad\text{and}\quad\phi^*R_P=A\,.$$
The Weyl conformal curvature tensor $W:=A-\sigma\rho$ is the projection of $A$ on $\ker(\rho)$; we say a model or a
pseudo Riemannian manifold is {\it conformally flat} if and only if $W=0$. The following
result \cite{BGKSW09} shows, in particular, that the relations of Equation (\ref{eqn-3.a}) generate the universal symmetries of the
Riemann curvature tensor. We focus our attention on the scalar curvature:

\begin{theorem}\label{thm-3.2}
\ \begin{enumerate}\item Any curvature model is geometrically realizable by a pseudo Riemannian manifold of constant
scalar curvature.
\item
Any conformally flat curvature model is geometrically realizable by a conformally flat pseudo Riemannian manifold of
constant scalar curvature.
\end{enumerate}\end{theorem}

To solve the {\it Yamabe problem} \cite{A76,R84,T68,Y60}, one constructs a Riemannian metric of constant scalar
curvature in a given conformal class on a compact Riemannian manifold. The complex analogue also has been solved \cite{RS03}
by constructing an almost Hermitian metric of constant scalar curvature in the conformal class of a compact almost Hermitian
manifold. Theorem \ref{thm-3.2} has a somewhat different flavor as we are not fixing the conformal class but rather the
curvature tensor at the realizing point. Furthermore, our manifolds are not compact nor even
complete.
\section{Affine and Riemannian Geometry}\label{S-4-AffineRiemann}
We now consider mixed structures -- we shall study an affine structure and a pseudo-Riemannian metric where the given affine
connection is not the Levi-Civita connection of the pseudo-Riemannian metric; thus the two structures are decoupled.

Let
$\langle\cdot,\cdot\rangle$ be a non-degenerate symmetric inner product on $V$ of signature $(p,q)$. Expand
$\mathcal{A}\in\mathfrak{A}(V)$ in the form
$$\mathcal{A}(e_i,e_j)e_k=\mathcal{A}_{ijk}{}^\ell e_\ell\,.$$
The scalar curvature
$\tau$ and trace free Ricci tensor $\rho_0$ are then given, respectively, by contracting indices
$$
\tau:=\varepsilon^{ij}\mathcal{A}_{kij}{}^k,\quad
\rho_0(x,y):=
\rho_s(x,y)-{\textstyle\frac{\tau}m}\langle x,y\rangle\,.
$$
 One has an orthogonal decomposition of $V^*\otimes V^*$ into irreducible
$O(V,\langle\cdot,\cdot\rangle)$ modules
\begin{equation}\label{eqn-4.a}
V^*\otimes V^*=\Lambda^2(V^*)\oplus S_0^2(V^*,\langle\cdot,\cdot\rangle)\oplus\mathbb{R}\,.
\end{equation}
\subsection{Geometric realizability of affine algebraic curvature tensors II} The decomposition of Equation (\ref{eqn-4.a})
leads to several geometric realization questions which are natural with respect to the structure group
$O(V,\langle\cdot,\cdot\rangle)$ and which can all be solved either in the real analytic category or in the $C^s$ category for any
$s\ge1$. As our considerations are local, we take $M=V$ and $P=0$. The primary focus of our investigation is on constant
scalar curvature and on properties of the Ricci tensor. We refer to \cite{GNW09p} for the proof of
the following result:

\begin{theorem}\label{thm-4.1}
Let $g$ be a $C^s$ (resp. real analytic) pseudo Riemannian metric which is defined on an open neighborhood of $0\in V$. Let
$\mathcal{A}\in\mathfrak{A}(V)$. There exists a torsion free
$C^s$ (resp. real analytic) connection $\nabla$ which is defined on a smaller open neighborhood of $0$ in $V$ such that:
\begin{enumerate}
\item $\mathcal{R}_0=\mathcal{A}$.
\item $\nabla$ has constant scalar curvature.
\item If $\mathcal{A}$ is Ricci symmetric, then $\nabla$ is Ricci symmetric.
\item If $\mathcal{A}$ is Ricci anti-symmetric, then $\nabla$ is Ricci anti-symmetric.
\item If $\mathcal{A}$ is Ricci traceless, then $\nabla$ is Ricci traceless.
\end{enumerate}\end{theorem}

\subsection{The decomposition of $\mathfrak{A}(V)$ as an $O(V,\langle\cdot,\cdot\rangle)$ module} The subspace
$\ker(\rho)\subset\mathfrak{A}(V)$ is not an irreducible
$O(V,\langle\cdot,\cdot\rangle)$ module but decomposes as the direct sum of 5 additional irreducible factors.
We refer to Bokan\cite{Bo90} for the proof of the following result:
\begin{theorem}\label{thm-4.2}
Let $\dim(V)\ge4$. There is an orthogonal decomposition of $\mathfrak{A}(V)$ into $8$ irreducible
$O(V,\langle\cdot,\cdot\rangle)$ modules
$\mathfrak{A}(V)= A_1\oplus...\oplus A_8$ where:
\medbreak\quad
$A_1\approx\mathbb{R}$,\quad
$A_2\approx A_5\approx S_0^2(V^*,\langle\cdot,\cdot\rangle)$,\quad
$A_3\approx A_4\approx\Lambda^2(V^*)$,
\smallbreak\quad
$A_6=\{\Theta\in\otimes^4V^*:\Theta_{ijkl}+\Theta_{jkil}+\Theta_{kijl}=0,
\Theta_{ijkl}=-\Theta_{jikl}=\Theta_{klij},$\smallbreak\qquad\qquad\qquad\qquad\quad$
\varepsilon^{il}\Theta_{ijkl}=0\}$,
\smallbreak\quad
$A_7=\{\Theta\in\otimes^4V^*:\Theta_{kjil}+\Theta_{ikjl}-\Theta_{ljik}-\Theta_{iljk}=0,
\Theta_{ijkl}=-\Theta_{jikl}=\Theta_{ijlk},$\smallbreak\qquad\qquad\qquad\qquad\quad$
  \textstyle\varepsilon^{il}\Theta_{ijkl}=0\}$,
\smallbreak\quad
$A_8=\{\Theta\in\otimes^4V^*:\Theta_{ijkl}=-\Theta_{jikl}=-\Theta_{ijlk}=-\Theta_{klij},
\varepsilon^{il}\Theta_{ijkl}=0\}$.
\end{theorem}

We note for the sake of completeness that:
$$\begin{array}{|l|l|}\gronk
\dim\{W_2\}=\dim\{W_5\}=\textstyle\frac{(m-1)(m+2)}2&
\dim\{W_1\}=1\\\gronk
\dim\{W_3\}=\dim\{W_4\}=\textstyle\frac{m(m-1)}2&
\dim\{W_6\}=\textstyle\frac{m(m+1)(m-3)(m+2)}{12}\\\gronk
\dim\{W_7\}=\textstyle\frac{(m-1)(m-2)(m+1)(m+4)}8&
\dim\{W_8\}=\textstyle\frac{m(m-1)(m-3)(m+2)}8\\\noalign{\hrule}
\end{array}$$

\section{almost Hermitian Geometry}\label{S5-almostHerm}
We refer to the discussion in \cite{DS88,F94,FFS,MC05,MC06,Sa03,S96} for additional information concerning almost Hermitian geometry.
We refer to \cite{CFG} for further information concerning almost para Hermitian geometry; for example,
para Hermitian geometry enters in the study of Osserman Walker metrics of signature $(2,2)$ \cite{nuevo}, it is important
in the study of homogeneous geometries
\cite{Gada}, and it is relevant to the study of Walker manifolds with degenerate self-dual Weyl curvature operators
\cite{Alex}.

Let $J$ be a
linear map of
$V$ and let
\
$\mathfrak{M}=(V,\langle\cdot,\cdot\rangle,A)$ be a curvature model.
One says that $J$ is a {\it Hermitian structure} if
$$J^2=-\operatorname{id}\quad\text{and}\quad J^*\langle\cdot,\cdot\rangle=\langle\cdot,\cdot\rangle\,.$$
Similarly, one says that $J$ is a {\it para Hermitian structure} if
$$J^2=\operatorname{id}\quad\text{and}\quad J^*\langle\cdot,\cdot\rangle=-\langle\cdot,\cdot\rangle\,.$$
Note that Hermitian structures exist if and only if both $p$ and $q$ are
even; para Hermitian structures exist if and only if $p=q$. Let
${\mathfrak{C}:=(V,\langle\cdot,\cdot\rangle,J,A)}$ be the associated {\it Hermitian curvature model} (resp. {\it
para Hermitian curvature model}). Define the $\star$-Ricci tensor $\rho^\star$ and the $\star$-scalar curvature $\tau^\star$
in either case, by setting:
\begin{eqnarray*}
&&\rho^\star(x,y):=\left\{
\begin{array}{lrl}
\phantom{-}\varepsilon^{il}A(e_i,x,Jy,Je_l)&\phantom{aaaa}\text{if}&\mathfrak{C}\text{ is Hermitian},\\
-\varepsilon^{il}A(e_i,x,Jy,Je_l)&\text{if}&\mathfrak{C}\text{ is para Hermitian},
\end{array}\right.\\
&&\tau^\star:=\phantom{Aaa.}\left\{
\begin{array}{lll}
\phantom{-}\varepsilon^{il}\varepsilon^{jk}A(e_i,e_j,Je_k,Je_l)&\text{if}&\mathfrak{C}\text{ is Hermitian},\\
-\varepsilon^{il}\varepsilon^{jk}A(e_i,e_j,Je_k,Je_l)&\text{if}&\mathfrak{C}\text{ is para Hermitian}\,.
\end{array}\right.
\end{eqnarray*}

\subsection{The geometric realizability of almost Hermitian models} One says that a manifold $\mathcal{C}:=(M,g,J)$ is an
{\it almost Hermitian manifold} (resp. {\it almost para Hermitian manifold}) if $(T_PM,g_P,J_P,R_P)$ is a
Hermitian (resp. para Hermitian) curvature model for every
$P\in M$. We do not assume that the structure $J$ on $M$ is integrable as this imposes additional curvature identities
\cite{gray} as we shall see presently. The notion of geometricrealizability in this context is defined similarly. Again, we focus our attention on the scalar curvature and the
analogous $\star$-scalar curvature in the following Theorem \cite{BGKSW09}:

\begin{theorem}\label{thm-5.1}
Let $m\ge4$.
\begin{enumerate}
\item Any Hermitian curvature model is geometrically realizable by an almost Hermitian manifold of
constant scalar and constant $\star$-scalar curvature.
\item Any para Hermitian curvature model is geometrically realizable by an
almost para Hermitian manifold of constant scalar and constant
$\star$-scalar curvature.
\end{enumerate}\end{theorem}

\subsection{The decomposition of $\mathfrak{R}(V)$ as a unitary module} Let
$(V,\langle\cdot,\cdot\rangle,J)$ be a Hermitian structure.
The Kaehler form is defined by
$\Omega(x,y):=\langle x,Jy\rangle$. Set
\medbreak\qquad\qquad
$S_{0,+}=S_{0,+}^2(V^*,\langle\cdot,\cdot\rangle,J):=\{\theta\in S^2(V^*):J^*\theta=\theta,\theta\perp\langle\cdot,\cdot\rangle\}$,
\smallbreak\qquad\qquad
$\Lambda_{0,+}=\Lambda_{0,+}^2(V^*,\langle\cdot,\cdot\rangle,J):=\{\theta\in\Lambda^2(V^*):J^*\theta=\theta,\theta\perp\Omega\}$,
\smallbreak\qquad\qquad
$S_-^2=S_-^2(V^*,J):=\{\theta\in S^2(V^*):J^*\theta=-\theta\}$,
\smallbreak\qquad\qquad
$\Lambda_-^2=\Lambda_-^2(V^*,J):=\{\theta\in\Lambda^2(V^*):J^*\theta=-\theta\}$\,.
\medbreak\noindent
\medbreak\noindent Let ${\mathcal{U}}(V,\langle\cdot,\cdot\rangle,J)$ be the associated unitary group:
$$\mathcal{U}(V,\langle\cdot,\cdot\rangle,J):=\{U\in\operatorname{GL}(V):UJ=JU\quad\text{and}\quad
U^*\langle\cdot,\cdot\rangle=\langle\cdot,\cdot\rangle\}\,.$$
We have the following decomposition of $V^*\otimes V^*$ as the orthogonal direct sum of
irreducible $\mathcal{U}(V,\langle\cdot,\cdot\rangle,J)$ modules:
$$V\otimes V^*=\langle\cdot,\cdot\rangle\cdot\mathbb{R}\oplus S_{0,+}^2\oplus
S_-^2\oplus\Omega\cdot\mathbb{R}\oplus\Lambda_{0,+}^2\oplus\Lambda_-^2.$$
We let  $\rho_{0,+,S}$, $\rho_{0,+,S}^\star$, $\rho_{-,S}$, and $\rho_{-,\Lambda}^\star$ denote the components of
$\rho$ and $\rho^\star$ with respect to this decomposition.

We refer to
\cite{TV81} for the proof of Theorem \ref{thm-5.2} in the Riemannian setting --  the extension to the higher signature
context is not difficult \cite{BGSV09}. This result  has been used by many authors
\cite{Bu07,F05,GM08,GM08a,MO08}.

\begin{theorem}\label{thm-5.2}
Let $(V,\langle\cdot,\cdot\rangle,J)$ be a Hermitian structure. \begin{enumerate}
\item We have the following orthogonal direct sum decomposition of $\mathfrak{R}(V)$ into irreducible
$\mathcal{U}(V,\langle\cdot,\cdot\rangle,J)$ modules:
\begin{enumerate}
\item If $2n=4$,
$\mathfrak{R}(V)=\WW_1\oplus\WW_2\oplus\WW_3\oplus\WW_4\oplus\WW_7
\oplus\WW_8\oplus\WW_9$.
\item If $2n=6$,
$\mathfrak{R}(V)=\WW_1\oplus\WW_2\oplus\WW_3\oplus\WW_4\oplus\WW_5\oplus\WW_7
\oplus\WW_8\oplus\WW_9\oplus\WW_{10}$.
\item If $2n\ge8$,
$\mathfrak{R}(V)=\WW_1\oplus\WW_2\oplus\WW_3\oplus\WW_4\oplus\WW_5\oplus\WW_6\oplus\WW_7
\oplus\WW_8\oplus\WW_9\oplus\WW_{10}$.
\end{enumerate}
We have $\WW_1\approx\WW_4$ and, if $2n\ge6$, $\WW_2\approx\WW_5$. The other $\mathcal{U}(V,\langle\cdot,\cdot\rangle,J)$
modules appear with multiplicity 1.
\item We have that:
\begin{enumerate}
\item  $\tau\oplus\tau^\star:\WW_1\oplus\WW_4\mapright{\approx}\mathbb{R}\oplus\mathbb{R}$.
\item If $2n=4$, $\rho_{0,+,S}:\WW_2\mapright{\approx} S_{0,+}^2$.
\item If $2n\ge6$, $\rho_{0,+,S}\oplus\rho_{0,+,S}^\star:\WW_2\oplus\WW_5\mapright{\approx} S_{0,+}^2\oplus S_{0,+}^2$.
\item $\WW_3
=\{A\in\mathfrak{R}(V):A(x,y,z,w)=A(Jx,Jy,z,w)\ \forall\ x,y,z,w\}\cap\ker(\rho)$.
\item If $2n\ge8$, $\WW_6=\ker(\rho\oplus\rho^\star)\cap\{A\in\mathfrak{R}(V):J^*A=A\}\cap\WW_3^\perp$.
\item $\WW_7=\{A\in\mathfrak{R}(V):A(Jx,y,z,w)=A(x,y,Jz,w)\ \forall\  x,y,z,w\}$.
\item $\rho_{-,S}:\WW_8\mapright{\approx} S_-^2$ and
$\rho_{-,\Lambda}^\star :\WW_9\mapright{\approx}\Lambda_-^2$.
\item If $2n\ge6$, $\WW_{10}=\{A\in\mathfrak{R}(V):J^*A=-A\}\cap\ker(\rho\oplus\rho^\star)$.
\end{enumerate}
\end{enumerate}
\end{theorem}

Let $m=\dim(V)=2n$. We note for the sake of completeness that:
\begin{equation}\label{eqn-5.a}
\begin{array}{|l|l|l|l||l|l|l|l|l|l|}\gronk&m=4&m=6&m\ge8&&m=4&m=6&m\ge8\\\gronk
\mathcal{W}_1&1&1&1&
\mathcal{W}_2&3&8&n^2-1\\\gronk
\mathcal{W}_3&5&27&\frac{n^2(n-1)(n+3)}4&
\mathcal{W}_4&1&1&1\\\gronk
\mathcal{W}_6&0&0&\frac{n^2(n+1)(n-3)}4&
\mathcal{W}_5&0&8&n^2-1\\\gronk
\mathcal{W}_7&2&12&\frac{n^2(n^2-1)}6&
\mathcal{W}_8&6&12&n^2+n\\\gronk
\mathcal{W}_{10}&0&30&\frac{2n^2(n^2-4)}3&
\mathcal{W}_9&2&6&n^2-n
\\\noalign{\hrule}
\end{array}\end{equation}

\subsection{The decomposition of $\mathfrak{R}(V)$ as a para unitary module}
We change the signs appropriately to obtain a corresponding decomposition in the para Hermitian context. Let
$(V,\langle\cdot,\cdot\rangle,\tilde J)$
 be a para Hermitian structure. Let $\tilde\Omega(x,y):=\langle
x,\tilde Jy\rangle$ be the para Kaehler  form. We have
$$\tilde J^*\tilde\Omega=-\tilde\Omega\quad\text{and}\quad
\tilde J^*\langle\cdot,\cdot\rangle=-\langle\cdot,\cdot\rangle\,.$$
Set\medbreak\qquad\qquad
$S_+^2=S_+^2( V^*,\tilde J):=\{\theta\in S^2( V^*):\tilde J^*\theta=\theta\}$,
\smallbreak\qquad\qquad
$\Lambda_+^2=\Lambda_+^2( V^*,\tilde J):=\{\theta\in\Lambda^2( V^*):\tilde J^*\theta=\theta\}$,
\smallbreak\qquad\qquad
$S_{0,-}^2=S_{0,-}^2( V^*,\langle\cdot,\cdot\rangle,\tilde J):=\{\theta\in S^2( V^*):\tilde
J^*\theta=-\theta,\theta\perp\langle\cdot,\cdot\rangle\}$,
\smallbreak\qquad\qquad
$\Lambda_{0,-}^2=\Lambda_{0,-}^2( V^*,\langle\cdot,\cdot\rangle,\tilde J):=\{\theta\in\Lambda^2( V^*):\tilde
J^*\theta=-\theta,\theta\perp\tilde\Omega\}$,
\smallbreak\qquad\qquad
$\tilde{\mathcal{U}}(V,\langle\cdot,\cdot\rangle,J):=\{\tilde U\in\operatorname{GL}( V):\tilde U\tilde
J=\tilde J\tilde U\quad\text{and}\quad
\tilde U^*\langle\cdot,\cdot\rangle=\langle\cdot,\cdot\rangle\}$.
\medbreak\noindent

\begin{theorem}\label{thm-5.3}
Let $(V,\langle\cdot,\cdot\rangle,\tilde J)$ be a para Hermitian structure. \begin{enumerate}
\item We have the following orthogonal direct sum decomposition of $\mathfrak{R}(V)$ into irreducible
$\tilde{\mathcal{U}}(V,\langle\cdot,\cdot\rangle,J)$ modules:
\begin{enumerate}
\item If $2n=4$,
$\mathfrak{R}(V)=\WWW_1\oplus\WWW_2\oplus\WWW_3\oplus\WWW_4\oplus\WWW_7
\oplus\WWW_8\oplus\WWW_9$.
\item If $2n=6$,
$\mathfrak{R}(V)=\WWW_1\oplus\WWW_2\oplus\WWW_3\oplus\WWW_4\oplus\WWW_5\oplus\WWW_7
\oplus\WWW_8\oplus\WWW_9\oplus\WWW_{10}$.
\item If $2n\ge8$,
$\mathfrak{R}(V)=\WWW_1\oplus\WWW_2\oplus\WWW_3\oplus\WWW_4\oplus\WWW_5\oplus\WWW_6\oplus\WWW_7
\oplus\WWW_8\oplus\WWW_9\oplus\WWW_{10}$.
\end{enumerate}
We have $\WWW_1\approx\WWW_4$ and, if $2n\ge6$, $\WWW_2\approx\WWW_5$. The other
$\tilde{\mathcal{U}}(V,\langle\cdot,\cdot\rangle,J)$ modules appear with multiplicity 1.
\item We have that:
\begin{enumerate}
\item  $\tau\oplus\tau^\star:\WWW_1\oplus\WWW_4\mapright{\approx}\mathbb{R}\oplus\mathbb{R}$.
\item If $2n=4$, $\rho_{0,-,S}:\WWW_2\mapright{\approx}S_{0,-}^2(V^*,\tilde{J})$.
\item If $2n\ge6$, $\rho_{0,-,S}\oplus\rho_{0,-,S}^\star:\WWW_2\oplus\WWW_5\mapright{\approx}S_{0,-}^2(V^*,\tilde{J})\oplus S_{0,-}^2(V^* ,\tilde{J})$.
\item $\WWW_3
=\{\tilde A\in\mathfrak{R}(V):\tilde A(x,y,z,w)=-\tilde A(\tilde Jx,\tilde Jy,z,w)\ \forall\
x,y,z,w\}$\newline$\phantom{AA}\cap\ker(\rho)$.
\item If $2n\ge8$, $\WWW_6=\ker(\rho\oplus\rho^\star)\cap\{\tilde A\in\mathfrak{R}(\tilde{V}):\tilde J^*\tilde A=\tilde
A\}\cap\WWW_3^\perp$.
\item $\WWW_7=\{\tilde A\in\mathfrak{R}(V):\tilde A(\tilde Jx,y,z,w)=\tilde A(x,y,\tilde Jz,w)\;\ \forall\
x,y,z,w\}$.
\item $\rho_{+,S}:\WWW_8\mapright{\approx}S_+^2(V^*,\tilde J)$,
$\rho_{+,\Lambda}^\star:\WWW_9\mapright{\approx}\Lambda_+^2(V^*,\tilde J)$.
\item If $2n\ge6$, $\WWW_{10}=\{\tilde A\in\mathfrak{R}(V):\tilde J^*\tilde A=-\tilde A\}\cap\ker(\rho\oplus\rho^\star)$.
\end{enumerate}
\end{enumerate}
\end{theorem}

We note for the sake of completeness that $\dim(\WWW_\nu)=\dim(\WW_\nu)$ is given in Equation (\ref{eqn-5.a}).
\section{Almost hyper Hermitian geometry}\label{S6-hyperHerm}
\subsection{Hyper Hermitian and hyper para Hermitian geometry} Fix a curvature model
$\mathfrak{M}=(V,\langle\cdot,\cdot\rangle,A)$. Let
$\mathcal{J}:=\{J_1,J_2,J_3\}$ be a triple of linear maps of $V$. We say that $\mathcal{J}$ is a {\it hyper Hermitian
structure} if $J_1$, $J_2$, $J_3$ are Hermitian structures and if we have the quaternion identities:
$$J_1^2=J_2^2=J_3^2=-\operatorname{id}\quad\text{and}\quad J_1J_2=-J_2J_1=J_3\,.$$
Similarly, we say that $\mathcal{J}$ is a {\it hyper para Hermitian structure} if $J_1$ is a Hermitian structure, if
$J_2$ and $J_3$ are para Hermitian structures, and if we have the para quaternion identities:
$$J_1^2=-J_2^2=-J_3^2=-\operatorname{id}\quad\text{and}\quad J_1J_2=-J_2J_1=J_3\,.$$
Let ${\mathfrak{Q}:=(V,\langle\cdot,\cdot\rangle,\mathcal{J},A)}$ be the associated
{\it hyper Hermitian curvature model} (resp. {\it hyper para Hermitian curvature model)}. We
refer to \cite{IZ2005,Kam99,CS04} for further details concerning such structures. We define:
$$\tau^\star_\mathfrak{Q}:=\tau^\star_{J_1}+\tau^\star_{J_2}+\tau^\star_{J_3}\,.$$

The structure group of a hyper Hermitian structure $\mathcal{J}$ is $SO(3)$ and of a hyper para Hermitian structure is
$SO(2,1)$ since we must allow for reparametrizations; $\tau^\star_\mathfrak{Q}$ is invariant under this structure group and does not
depend on the particular parametrization chosen. We say that $(M,g,\mathcal{J})$ is an  {\it almost
hyper Hermitian} manifold  or an  {\it almost hyper para Hermitian manifold} if $\mathcal{J}_P$ defines the
appropriate structure on $(T_PM,g_P)$ for all points
$P$ of
$M$; we impose no integrability condition.

\begin{theorem}\label{thm-6.1}
Let $m\ge8$.
\begin{enumerate}
\item Any hyper Hermitian curvature model is geometrically realizable by an almost hyper Hermitian
manifold of constant scalar and constant
$\star$-scalar curvature.
\item Any hyper para Hermitian curvature model is geometrically realizable by an almost hyper para
Hermitian manifold of constant scalar and constant
$\star$-scalar curvature.
\end{enumerate}
\end{theorem}

\section{Hermitian Geometry}\label{S7-HermGeo}
We refer to
\cite{AGI97, B90,JK07,Sa04,V07} for additional material on Hermitian geometry.
We say an almost Hermitian manifold $\mathcal{M}=(M,g,\mathcal{J})$ is {\it Hermitian} if $\mathcal{J}$ is an
integrable almost complex structure, i.e. the {\it Nijenhuis tensor}
$$N_{\mathcal{J}}(x,y):=[x,y]+\mathcal{J}[\mathcal{J}x,y]+\mathcal{J}[x,\mathcal{J}y]-[\mathcal{J}x,\mathcal{J}y]$$
 vanishes
or, equivalently, in a neighborhood of any point of the manifold there are local coordinates
$(x_1,...,x_n,y_1,...,y_n)$ so that
$$\mathcal{J}\partial_{x_i}=\partial_{y_i}\quad\text{and}\quad \mathcal{J}\partial_{y_i}=-\partial_{x_i}\,.$$
Similarly \cite{CMMS03}, we say that $(M,g,\tilde{\mathcal{J}})$ is a {\it para Hermitian manifold} if
$\tilde{\mathcal{J}}$ is an integrable almost para complex structure, i.e.  if the {\it para Nijenhuis tensor}
$N_{\tilde{\mathcal{J}}}$
$$N_{\tilde{\mathcal{J}}}(x,y):=[x,y]-\tilde{\mathcal{J}}[\tilde{\mathcal{J}}x,y]
-\tilde{\mathcal{J}}[x,\tilde{\mathcal{J}}y]+[\tilde{\mathcal{J}}x,\tilde{\mathcal{J}}y]$$
 vanishes
or, equivalently, there exist local coordinates $(x_1,...,x_n,y_1,...,y_n)$
centered at any given point of $M$ so  that
$$\tilde{\mathcal{J}}\partial_{x_i}=\partial_{y_i}\quad\text{and}\quad\tilde{\mathcal{J}}\partial_{y_i}=\partial_{x_i}\,.$$

Gray
\cite{gray} showed that the curvature tensor of a Hermitian manifold has an additional symmetry given below in Equation (\ref{eqn-7.a});
it is quite striking that a geometric integrability condition imposes an additional algebraic symmetry on the curvature
tensor. We refer to \cite{gray} for the proof of the first Assertion and to \cite{BGSV09} for the proof of the second Assertion
in the following Theorem:

\begin{theorem}\label{thm-7.1}
\ \begin{enumerate}
\item If a  Hermitian curvature model $\mathfrak{C}=(V,\langle\cdot,\cdot\rangle,J,A)$ is geometrically realizable by a
Hermitian manifold, then
\begin{eqnarray}
0&=&A(x,y,z,w)+A(Jx,Jy,Jz,Jw)\nonumber\\
&+&A(Jx,Jy,z,w)+A(x,y,Jz,Jw)+A(Jx,y,Jz,w)\label{eqn-7.a}\\
&+&A(x,Jy,z,Jw)+A(Jx,y,z,Jw)+A(x,Jy,Jz,w)\}\,.\nonumber
\end{eqnarray}
\item If a para Hermitian curvature model $\tilde{\mathfrak{C}}=(V,\langle\cdot,\cdot\rangle,\tilde J,\tilde A)$ is
geometrically realizable by a para Hermitian manifold, then
\begin{eqnarray}
0&=&\tilde A(x,y,z,w)+\tilde A(\tilde Jx,\tilde Jy,\tilde Jz,\tilde Jw)\nonumber\\
&-&\tilde A(\tilde Jx,\tilde Jy,z,w)-\tilde A(x,y,\tilde Jz,\tilde Jw)-\tilde A(\tilde Jx,y,\tilde Jz,w)\label{eqn-7.b}\\
&-&\tilde A(x,\tilde Jy,z,\tilde Jw)-\tilde A(\tilde Jx,y,z,\tilde Jw)-\tilde A(x,\tilde Jy,\tilde Jz,w)\}\,.\nonumber
\end{eqnarray}
\end{enumerate}\end{theorem}

We refer to \cite{BGKW09} for the proof of the first Assertion and to \cite{BGSV09} for the proof of the second Assertion in
Theorem \ref{thm-7.2}; this result provides a useful converse to Theorem \ref{thm-7.1}. Again we shall focus on the scalar
curvature and the
$\star$-scalar curvature:
\begin{theorem}\label{thm-7.2}
\ \begin{enumerate}
\item If a Hermitian curvature model $\mathfrak{C}$ satisfies {\rm Equation (\ref{eqn-7.a})}, then $\mathfrak{C}$ is
geometrically realizable by a Hermitian manifold with constant scalar curvature, with constant
$\star$-scalar curvature, and with $d\Omega$ vanishing at the realizing point $P$.
\item If a para Hermitian curvature model $\tilde{\mathfrak{C}}$ satisfies {\rm Equation (\ref{eqn-7.b})}, then
$\tilde{\mathfrak{C}}$ is geometrically realizable by a para Hermitian manifold with constant scalar
curvature, with constant $\star$-scalar curvature, and with
$d\tilde\Omega$ vanishing at the realizing point $P$.
\end{enumerate}\end{theorem}

Equation (\ref{eqn-7.a}) is called the {\it Gray identity} and Equation
(\ref{eqn-7.b}) is called the {\it para Gray identity}. The
universal symmetries of the curvature tensor of a Hermitian manifold (resp. a para Hermitian
manifold) are generated by the Gray (resp. para Gray) identity and the usual curvature symmetries
(see Equation (\ref{eqn-3.a})). This result emphasizes the difference between almost Hermitian and
Hermitian manifolds.

\begin{remark}\label{rmk-7.3}
\rm Since the Hermitian geometric
realization can be chosen so that $d\Omega(P)=0$, imposing the Kaehler identity $d\Omega (P)=0$ at a
single point imposes no additional curvature restrictions. If $d\Omega=0$ globally, then the manifold is said
to be almost Kaehler. This is a very rigid structure, see for example the discussion in
\cite{T06}, and there are additional curvature restrictions. Thus Theorem \ref{thm-7.2} also emphasizes the difference
between
$d\Omega$ vanishing at a single point and
$d\Omega$ vanishing globally.\end{remark}

\begin{remark}\label{rmk-7.4}\rm
The space of vectors
satisfying the Gray condition (resp. the para Gray condition) is exactly
$\mathcal{W}_7^\perp$ (resp. $\WWW_7^\perp$) of Theorem \ref{thm-5.2} (resp. Theorem \ref{thm-5.3}). Furthermore, either the complex
Jacobi operator or the complex curvature operator completely determine the components in
$\mathcal{W}_7^\perp$ of a curvature tensor \cite{BVGGR08}; the algebraic condition
determining $\mathcal{W}_7$ also plays a role in the study of Jacobi--Ricci commuting curvature tensors \cite{GN}.
\end{remark}

\section{Ivanov--Petrova geometry}\label{sect-9-IvanovPetrovaGeometry} To simplify the discussion, we work in the Riemannian setting;
there are analogous results in arbitrary signatures. Let
$\mathfrak{M}=(V,\langle\cdot,\cdot\rangle,A)$ be a Riemannian curvature model.
Let $\{x,y\}$ be a basis for an oriented $2$-plane
$\pi$. The {\it skew-symmetric curvature operator} $\mathcal{R}(\pi)$ is defined by setting
$$\mathcal{R}(\pi):=\{\langle x,x\rangle\langle y,y\rangle-\langle x,y\rangle^2\}^{-1/2}\mathcal{R}(x,y)\,.$$
This skew-symmetric operator is independent of the particular basis chosen. We say that $\mathfrak{M}$ isIvanov--Petrova if the eigenvalues of
$\mathcal{R}(\pi)$ are constant on the Grassmannian of oriented $2$-planes in $V$. Similarly we say that a Riemannian manifold
$\mathcal{M}=(M,g)$ is Ivanov--Petrova if the curvature model $(T_PM,g_P,R_P)$  is Ivanov--Petrova for all points $P$ of the manifold.

The study of such manifolds was initiated by Ivanov and Petrova \cite{IP98} in dimension $4$ and the notation ``Ivanov--Petrova'' was
adopted by later authors following this seminal paper. Let $\phi$ be a self-adjoint map, i.e. $\langle\phi x,y\rangle=\langle x,\phi
y\rangle$ for all
$x,y$. We form the following algebraic curvature tensor:
\begin{equation}\label{eqn-8.x}
A_\phi(x,y,z,w)=\{\langle \phi x,z\rangle\langle\phi y,zw\rangle-\langle\phi x,w\rangle\langle\phi y,z\rangle\}\,.
\end{equation}

\subsection{Ivanov-Petrova curvature models} We have the following examples \cite{G-99,GLS99,IP98} in the algebraic context:

\begin{example}\label{exm-8.1}
\rm\ \begin{enumerate}
\item Let $\phi$ be a self-adjoint map of $(V,\langle\cdot,\cdot\rangle)$ with $\phi^2=\operatorname{id}$.
Adopt the notation of Equation (\ref{eqn-8.x}). Then $(V,\langle\cdot,\cdot\rangle,CA_\phi)$ is Ivanov--Petrova for any constant $C$.
Note that if $\phi=\pm\operatorname{id}$, then $A_\phi$ has constant sectional curvature $C$.
\item Let $\{e_1,e_2,e_3,e_4\}$ be the standard normalized orthonormal basis for
$\mathbb{R}^4$. Let $2a_1+a_2=0$. Let the non-zero components of $A$, up to the
usual symmetries, be:
$$\begin{array}{rrrr}
     R_{1212}=a_1,& R_{1234}=\phantom{-}a_2,&
     R_{1313}=a_2,& R_{1324}=-a_1,\\
     R_{1414}=a_2,& R_{1423}=\phantom{-}a_1,&
     R_{2323}=a_2,& R_{2314}=\phantom{-}a_1,\\
     R_{2424}=a_2,& R_{2413}=-a_1,&
     R_{3434}=a_1,\hfill&R_{3412}=\phantom{-}a_2\,.\end{array}$$
This tensor is Ivanov--Petrova; we also refer to
\cite{G-99} where this tensor is described in terms of quaternions.
\end{enumerate}\end{example}

 One has a complete classification. The 4 dimensional case is exceptional and is covered by 
Ivanov and Petrova \cite{IP98}. The
cases
$m\ge5$ and
$m\ne7,8$ are dealt with by the work of \cite{GLS99}. The case $m=8$ is treated in \cite{G-99} and the case $m=7$ is
presented in Nikolayevsky
\cite{N04}.

\begin{theorem}\label{thm-8.2}
Any Ivanov--Petrova curvature model of dimension $m\ge4$ is
isomorphic to one of the models of {\rm Example \ref{exm-8.1}}.
\end{theorem}

\subsection{Ivanov-Petrova manifolds} When we pass to the geometric setting, we have the following examples
\cite{GLS99,IP98} of Ivanov--Petrova manifolds:
\begin{example}\label{exm-8.3}
\rm\ \begin{enumerate}
\item Any manifold of constant sectional curvature is Ivanov--Petrova.
\item Let $(N,h)$ be a metric of constant sectional curvature $K\ne0$. Consider $f(t)=Kt^2+At+B$ and let $\mathcal{O}$ be an
open subset of $\mathbb{R}$ where $f>0$. Let
$M=\mathcal{O}\times N$ with the metric
$$g=dt^2+f(t)h\,.$$
If $A^2-4KB\ne0$, this metric does not have constant sectional curvature and is Ivanov-Petrova.
\end{enumerate}
\end{example}

We then have the following geometric classification Theorem \cite{G-99,GLS99,IP98,N04}:

\begin{theorem}\label{thm-8.4}
Any Ivanov--Petrova manifold of dimension $m\ge4$ is locally isometric to one of
the manifolds of {\rm Example \ref{exm-8.3}}.
\end{theorem}

\subsection{Ivanov-Petrova curvature models which are not geometrically realizable by Ivanov-Petrova manifolds} If $\phi$ is
a self-adjoint map of
$(V,\langle\cdot,\cdot\rangle)$ with
$\phi^2=\operatorname{id}$, we can find an orthonormal basis
$\{e_i\}$ for
$V$ so that
$\phi e_i=\pm e_i$. Let $p$ be the number of $+1$ eigenvalues and $q$ the number of $-1$ eigenvalues; the modified inner product
$\langle x,y\rangle_\phi:=\langle\phi x,y\rangle$ has signature $(p,q)$. The curvature tensors of the manifolds of Example \ref{exm-8.3}
(1) correspond to $\phi$ with $\langle\cdot,\cdot\rangle_\phi$ having signature $(0,m)$ or $(m,0)$ (i.e.
$\phi=\pm\operatorname{id}$) and the curvature tensors of Example
\ref{exm-8.3} (2) correspond to $\phi$ with $\langle\cdot,\cdot\rangle_\phi$ having signature $(1,m-1)$ or $(m-1,1)$.
Thus we may combine Theorems \ref{thm-8.2} and Theorem
\ref{thm-8.4} to see:

\begin{theorem}\label{thm-8.6}
Any Ivanov--Petrova curvature model of dimension $m\ge4$ is geometrically realizable by an
Ivanov--Petrova manifold if and only if it has the form given in {\rm Example \ref{exm-8.1} (1)} where
$\langle\cdot,\cdot\rangle_\phi$ has signature
$(0,m)$,
$(m,0)$, $(1,m-1)$, or $(m-1,1)$.
\end{theorem}

Note that the exceptional Ivanov-Petrova model of Example \ref{exm-8.1} (2) in dimension $4$ is not geometrically
realizable by an Ivanov--Petrova manifold.

\section{Osserman Geometry}\label{sect-10-OssermanGeometry}
Fix a curvature model $\mathfrak{M}=(V,\langle\cdot,\cdot\rangle,A)$; again, we restrict to the Riemannian setting.
Let $S$ be the unit
sphere in
$(V,\langle\cdot,\cdot\rangle)$. The Jacobi operator
$$\mathfrak{J}(x):y\rightarrow\mathcal{R}(y,x)x$$
is a self-adjoint operator which
appears in the study of geodesic sprays. We say that a Riemannian curvature model or a Riemannian manifold is {\it Osserman}
if the eigenvalues of
$\mathfrak{J}$ are constant on $S$.

\subsection{Osserman curvature models} Let $\psi$ be a skew-adjoint map of $(V,\langle\cdot,\cdot\rangle)$. Motivated by
the splitting $\sigma\rho_a$ of Theorem \ref{thm-2.1}, we form the algebraic curvature tensor:
$$A_\psi(x,y,z,w):=\langle\psi y,z)\langle\psi x,w\rangle-\langle\psi x,z\rangle\langle\psi y,w\rangle
-2\langle\psi x,y\rangle\langle\psi z,w\rangle\,.$$
The following examples appear first in \cite{G94}:

\begin{example}\label{ex-9.1}
\rm Let $\{\psi_1,...,\psi_\ell\}$ be a family of skew-adjoint endomorphisms
defined on
$(V,\langle\cdot,\cdot\rangle)$ which satisfy the {\it Clifford commutation
relations}:
$$\psi_i\psi_j+\psi_j\psi_i=-2\delta_{ij}\,.$$
Let $\{\lambda_0,\lambda_1,...,\lambda_\ell\}$ be real constants where $\lambda_i\ne0$ if $i>0$. Set
$$A=\lambda_0A_{\operatorname{id}}+\lambda_1A_{\psi_1}+...+\lambda_\ell A_{\psi_\ell}\,.$$
Then $(V,\langle\cdot,\cdot\rangle,A)$ is an Osserman curvature model. If
$\ell=0$, then
$A=\lambda_0A_{\operatorname{id}}$ has constant sectional curvature.
\end{example}

\begin{remark}\label{rmk-9.2}
\rm The family $\{\psi_i\}$ is said to give a {\it Clifford module} structure to $V$ and the maximal such $\ell$ possible is
called the {\it Adams number} and is denoted by $\nu(m)$. If $m$ is odd, no such structure is possible and $\nu(m)=0$.
If
$m\equiv2$ mod $4$, then only $\ell=1$ is possible. If $m\equiv4$ mod $8$, then $\ell=3$ is possible; this case can be
realized by a quaternion structure (although there are other structures possible if $m>4$). We refer to Adams \cite{A62} for
further details as this number is closely related to the number of linearly independent vector fields on spheres. If
$m=a2^s$ where
$a$ is odd, then
$\nu(m)=\nu(2^s)$. We have
$$\nu(1)=0,\quad\nu(2)=1,\quad\nu(4)=3,\quad\nu(8)=7,\quad\nu(16\cdot 2^k)=8+\nu(2^k)\,.$$
\end{remark}

\begin{theorem}\label{thm-9.3}
Let $\mathfrak{M}$ be an Osserman curvature model of dimension $m\ne16$. Then $\mathfrak{M}$ is isomorphic to one of the
curvature models of {\rm Example \ref{ex-9.1}}.
\end{theorem}

Theorem \ref{thm-9.3} was proved by Chi \cite{C88} if $m\equiv1$ mod $2$, if $m\equiv 2$ mod $4$, and if $m=4$. Subsequently
Nikolayevsky
\cite{N03,Ni1,Ni2} established Theorem \ref{thm-9.3} for the remaining values. The result fails in dimension 16; the
curvature model of the Cayley plane is not given by a Clifford module
structure and the classification is unknown in that dimension.

\subsection{Osserman manifolds} Osserman \cite{O90} conjectured that any Riemannian manifold whose Jacobi operator had
constant eigenvalues on the set of unit tangent vectors was necessarily a local $2$-point homogeneous space, i.e. is either
flat or is a rank 1 symmetric space. This conjecture became known as the Osserman conjecture by subsequent authors and thecondition that the Jacobi operator has constant eigenvalues is known as the Osserman condition. The Osserman conjecture in
the Riemannian setting has been settled except in dimension 16 where it remains open
\cite{C88,N03,Ni1,Ni2}:
\begin{theorem}\label{thm-9.4}
If $\mathcal{M}$ is a Riemannian Osserman manifold of dimension $m\ne16$, then either $\mathcal{M}$ is flat or $\mathcal{M}$ is locally
isometric to a rank $1$ symmetric space.
\end{theorem}

\begin{remark}\label{rmk-9.5}
\rm Theorem \ref{thm-9.4} fails in the indefinite setting. There are Walker manifolds of signature $(2,2)$ which are Osserman
but not locally  homogeneous \cite{BGGNV,GKV}.\end{remark}

\subsection{Osserman curvature models which are not geometrically realizable by Osserman manifolds} We adopt the notation of
Example
\ref{ex-9.1}. If
$\mathfrak{M}$ is the curvature model of a rank
$1$ symmetric space ot dimension
$m\ne16$, then one of the following cases holds:
\begin{enumerate}\item
$\ell=0$ and
$\mathcal{M}$ has constant sectional curvature so $\mathcal{M}$ is locally isometric to a rescaled sphere or hyperbolic space.
\item $\ell=1$
and
$\lambda_1=\lambda_0$ so
$\mathcal{M}$ is locally isometric to a rescaled complex projective space or the negative curvature dual.
\item
$\ell=3$ and
$\lambda_0=\lambda_1=\lambda_2=\lambda_3$ so $\mathcal{M}$ is locally isometric to a rescaled quaternionic projective space or the
negative curvature dual.
\end{enumerate}

The following result now follows from the discussion above:

\begin{theorem}\label{thm-9.6}
Let $\mathfrak{M}$ be an Osserman curvature model of dimension $m\ne16$. Adopt the notation of
{\rm Example \ref{ex-9.1}}. Assume $\ell\ge1$. If
$\ell\ne1,3$ or if $\lambda_i\ne\lambda_j$ for some $i\ne j$, then $\mathfrak{M}$ is not geometrically realizable by an Osserman manifold.
\end{theorem}

\section{Curvature homogeneity}\label{sect-11-Curvaturehomogeneity}

Let $\mathfrak{M}=(V,\langle\cdot,\cdot\rangle)$ be a curvature model. We say that a pseudo Riemannian manifold
$\mathcal{M}=(M,g)$ is {\it curvature homogeneous} with model $\mathfrak{M}$ if every point of $M$ realizes $\mathfrak{M}$
geometrically; in this situation,
$\mathcal{M}$ is said to be {\it curvature homogeneous}. Equivalently, this means that given any $2$ points $P$ and $Q$ of $M$,
there is an isometry $\Phi_{P,Q}:T_PM\rightarrow T_QM$ so $\Phi_{P,Q}^*R_Q=R_P$. More generally, if $\Phi_{P,Q}^*\nabla^iR_Q=\nabla^iR_P$
for $i\le k$, then $\mathcal{M}$ is said to be {\it $k$ curvature homogeneous}. One has the following result of Singer
\cite{S60} in the Riemannian setting and of Podesta and Spiro \cite{PS96} in the pseudo Riemannian setting:
\begin{theorem}\label{thm-10.1}
 There exists an integer $k_{p,q}$ so that if $\mathcal{M}$ is any complete simply connected
pseudo Riemannian manifold of signature $(p,q)$ which is $k_{p,q}$ curvature homogeneous, then $\mathcal{M}$ is homogeneous.
\end{theorem}

One has rigidity results of Tricerri and Vanhecke \cite{TV86} and of Cahen, Leroy, Parker, Tricerri, and
Vanhecke \cite{CLPTV91}:

\begin{theorem}\label{thm-10.2}
\ \begin{enumerate}
\smallbreak\item
A Riemannian curvature homogeneous manifold which is $0$ curvature modeled on an
irreducible symmetric space is locally symmetric.
\smallbreak\item A Lorentzian curvature homogeneous
manifold which is $0$ curvature modeled on an irreducible symmetric space has constant sectional
curvature.\end{enumerate}\end{theorem}

In the Riemannian setting, \cite{FKM81,Ty74} there are curvature
homogeneous manifolds which are not locally homogeneous but there are no known examples which are
$1$ curvature homogeneous but not locally homogeneous. Work of \cite{S95,SSVa92} shows that any $1$ curvature
homogeneous complete simply connected Riemannian manifold of dimension $m\le5$ is homogeneous.
In the Lorentzian setting ($p=1$), there exist $1$ curvature homogeneous Lorentzian manifolds which are not locally
homogeneous
\cite{BD00,BV97}. On the other hand, given any $k$, one can construct neutral signature pseudo Riemannian manifolds which
are complete, which are modeled on a symmetric space, which are $k$ curvature homogeneous, and which are
not locally homogeneous \cite{GNc04d}.

The results of Sections \ref{sect-9-IvanovPetrovaGeometry} and \ref{sect-10-OssermanGeometry} immediately yield:
\begin{theorem}\label{thm-10.3}
There exist Riemannian curvature models which are not geometrically realizable by curvature homogeneous manifolds.
\end{theorem}

\section*{Acknowledgments}
Research of M. Brozos-V\'azquez was partially supported by Project MTM2006-01432 (Spain).
Research of P. Gilkey partially supported Projects MTM2006-01432 (Spain) and DGI SEJ2007-67810a (Spain).
Research of S. Nik\v cevi\'c partially supported by Research of Project 144032 (Serbia).


\begin{thebibliography}{AAA}

\bibitem{A62} J. Adams, {\it Vector fields on spheres}, Annals of Math. {\bf 75} (1962), 603--632.

\bibitem{AGI97} V. Apostolov, G. Ganchev, and S. Ivanov, {\it Compact Hermitian surfaces of
constant antiholomorphic sectional curvatures}, {Proc. Amer. Math. Soc.  \bf 125}  (1997), 3705--3714.


\bibitem{A76}T. Aubin, {\it \'Equations diff\'erentielles non lin\'eaires et probl\'eme de Yamabe concernant la courbure
scalaire},  {J. Math. Pures Appl.} {\bf 55} (1976), 269--296.


\bibitem{B90} D. E. Blair, {\it Nonexistence of $4$-dimensional almost Kaehler manifolds of constant curvature},
{Proc. Amer. Math. Soc.} {\bf 110} (1990), 1033--1039.

\bibitem{BGNS06} N.  Bla{\v z}i{\'c},  P. Gilkey, S. Nik\v cevi\'c, and U. Simon,
{\it Algebraic theory of affine curvature tensors},
Archivum Mathematicum (BRNO), \textbf{42} (2006),
Suppl., 147--168.

\bibitem{Bo90} N. Bokan,
{\it On the complete decomposition of curvature tensors of Riemannian manifolds with symmetric connection},Rend. Circ. Mat. Palermo \textbf{XXIX} (1990), 331--380.

\bibitem{BDS03} N. Bokan, M. Djori\'c, and U. Simon,
{\it Geometric structures as determined by the volume of generalized geodesic balls},
Result. Math. \textbf{43} (2003), 205--234.

\bibitem{BVGGR08} M. Brozos-V\'azquez, E. Garc\'ia-R\'io, and P. Gilkey,
{\it Relating the curvature tensor and the complex Jacobi operator of an almost Hermitian manifold}, {Adv. Geom.}
{\bf 8}  (2008), 353--365.

\bibitem{BGGNV} M. Brozos-V\'{a}zquez, E. Garc\'{\i}a-R\'{\i}o,
P. Gilkey, S. Nik\v cevi\'c, and R. V\'{a}zquez-Lorenzo
{\it The Geometry of Walker Manifolds}, Morgan and Claypool (to appear).

\bibitem{BGKW09} M.  Brozos-V\'azquez, P. Gilkey, H. Kang, and S. Nik\v cevi\'c, {\it Geometric realizations of Hermitian
curvature models}, arXiv:0812.2743.

\bibitem{BGKSW09} M.  Brozos-V\'azquez, P. Gilkey, H. Kang, S. Nik\v cevi\'c, and G. Weingart,
{\it Geometric realizations of curvature models by manifolds with constant scalar curvature}, 
to appear Differential Geom. Appl., arXiv:0811.1651.

\bibitem{BGSV09} M.  Brozos-V\'azquez, P. Gilkey, S. Nik\v cevi\'c, and  R. V\'{a}zquez-Lorenzo, {\it Geometric Realizations of
para Hermitian curvature models}, arXiv:0902.1697.

\bibitem{BD00} P. Bueken, and M. Djoric,
{\it Three-dimensional Lorentz metrics and curvature homogeneity of order one}, {Ann. Global
Anal. Geom.} {\bf 18} (2000), 85--103.

\bibitem{BV97} P. Bueken, and L. Vanhecke,
{\it Examples of curvature homogeneous Lorentz metrics},
{Classical Quantum Gravity} {\bf 14} (1997),
L93--L96.

\bibitem{Bu07} J. B. Butruille, {\it Espace de twisteurs d une vari\'et\'e presque hermitienne de dimension 6},
{Ann. Inst. Fourier (Grenoble) \bf 57} (2007), 1451--1485.

\bibitem{CLPTV91}
M. Cahen, J. Leroy, M. Parker, F. Tricerri, and L. Vanhecke,
{\it Lorentz manifolds modeled on a Lorentz symmetric space},
{J. Geom. Phys.} {\bf 7} (1990), 571--581.


\bibitem{C88} Q. S. Chi,
{\it A curvature characterization of certain locally rank one symmetric spaces},
J. Differential Geom. \textbf{28}
(1988), 187--202.


\bibitem{Alex}
 A. Cort\'{e}s-Ayaso,  J. C. D\'{\i}az-Ramos, and E. Garc\'{\i}a-R\'{\i}o, {\it Four-dimensional manifolds
with non-degenerated self-dual Weyl curvature tensor}, {Ann.
Global Anal. Geom.} {\bf 34} (2008), 185--193.

\bibitem{CMMS03} V. Cort\'es, C. Mayer, T. Mohaupt, and F. Saueressig, {\it Special geometry of Euclidean supersymmetry I: vector
multiplets}, arXiv:hep-th/0312001.

\bibitem{CFG}
V. Cruceanu, P. Fortuny, and P. M. Gadea, {\it A survey on paracomplex
geometry}, {Rocky Mount. J. Math.} {\bf 26} (1996), 83--115.

\bibitem{RS03} {H. del Rio, and S. Simanca, {\it The Yamabe problem for almost Hermitian manifolds},
{J. Geom. Anal.} {\bf 13} (2003), 185--203.}

\bibitem{DS88} J. Deprez, K. Sekigawa, and L. Verstraelen, {\it Classifications of Kaehler manifolds satisfying some
curvature conditions}, Sci. Rep. Niigata Univ. Ser. A {\bf 24} (1988), 1--12.

\bibitem{nuevo}
J. C. D\'{\i}az-Ramos, E. Garc\'{\i}a-R\'{\i}o, and R. V\'{a}zquez-Lorenzo,
{\it Osserman metrics on Walker $4$-manifolds equipped with a
para Hermitian structure},
{Mat. Contemp.} {\bf 30} (2006), 91--108.

\bibitem{F94} M. Falcitelli, and A. Farolina {\it Curvature properties of almost Hermitian manifolds}, Riv. Mat.  Univ.
Parma V (1994), 301--320.

\bibitem{FFS} M. Falcitelli, A. Farinola, and S. Salamon, {\it Almost-Hermitian geometry}, {Differential Geom. Appl.}
{\bf  4} (1994),  259--282.

\bibitem{FKM81}
D. Ferus, H. Karcher, and H. M\"unzner,
{\it Cliffordalgebren und neue isoparametrische Hyperfl\"achen},
{Math. Z.} {\bf 177} (1981), 479--502.

\bibitem{F05}A. Fino, {\it Almost Kaehler 4-dimensional Lie groups with
J-invariant Ricci tensor}, {Differential Geom. Appl. \bf 23} (2005), 26--27.

\bibitem{Gada}
P. M. Gadea, and J. A. Oubi\~{n}a, {\it Homogeneous pseudo-Riemannian structures
and homogeneous almost para Hermitian structures}, {Houston J.
Math.} {\bf 18} (1992),  449--465.

\bibitem{GM08} G. Ganchev, and V. Mihova,
{\it Kaehler manifolds of quasi-constant holomorphic sectional curvatures},
{Cent. Eur. J. Math. \bf 6} (2008), 43--75.

\bibitem{GM08a} G. Ganchev, and V. Mihova, {\it Warped product Kaehler manifolds and Bochner-Kaehler metrics},
{J. Geom. Phys. \bf 58} (2008), 803--824.

\bibitem{GKV}
E. Garc\'{\i}a-R\'{\i}o, D. Kupeli, and  R. V\'{a}zquez-Lorenzo,
{\it Osserman manifolds in semi-Riemannian geometry}, Lect. Notes
Math. {\bf 1777}, {Springer-Verlag}, Berlin,
2002.

\bibitem{G94}
P. Gilkey, {\it Manifolds whose curvature operator has
constant eigenvalues at the basepoint}, {J. Geom. Anal.} {\bf
4} (1994), 155--158.

\bibitem{G-99} P. Gilkey, {\it Riemannian manifolds whose skew-symmetric curvature operator has constant eigenvalues II},
Differential Geometry and applications  (eds Kolar, Kowalski, Krupka, and Slovak) Publ Masaryk University Brno Czech Republic ISBN
80-210-2097-0 (1999), 73--87.

\bibitem{GLS99} P. Gilkey, J. V. Leahy, and H. Sadofsky,
{\it Riemannian manifolds whose skew-symmetric curvature operator as constant
eigenvalues}, Indiana Univ. Math. J. {\bf 48} (1999), 615--634.

\bibitem{GNc04d}
P. Gilkey, and S. Nik\v cevi\'c,
{\it Complete $k$-curvature homogeneous pseudo-Riemannian manifolds},
{Ann. Global Anal. Geom.} {\bf 27} (2005), 87--100.


\bibitem{GN} P. Gilkey, and S. Nik\v cevi\'c, {\it Pseudo-Riemannian Jacobi--Videv Manifolds}, {
 Int. J. Geom. Methods Mod. Phys.} {\bf  4} (2007),  727--738.

\bibitem{GN08} P. Gilkey, and S. Nik\v cevi\'c, {\it Geometrical representations of equiaffine curvature operators},
Result. Math. {\bf 52} (2008), 281--287.



\bibitem{GNS09} P. Gilkey, S. Nik\v cevi\'c, and U. Simon, {\it Geometric theory of equiaffine curvature tensors},
arXiv:0903.5269.


\bibitem{GNW09}  P. Gilkey, S. Nik\v cevi\'c, and D. Westerman, {\it Geometric realizations of generalized algebraic curvature
operators}, 	J. Math. Phys. {\bf 50}, 013515 (2009).

\bibitem{GNW09p} P. Gilkey, S. Nik\v cevi\'c, and D. Westerman, {\it Riemannian geometric realizations for Ricci tensors
of generalized algebraic curvature operators}, to appear in the Conference Proceedings of the VIII International Colloquium on
Differential Geometry (World Sci. Publ. Co.).

\bibitem{gray} A. Gray, {\it Curvature identities for Hermitian and
almost Hermitian manifolds}, {T{\^o}hoku Math. J.} {\bf 28}
(1976), 601--612.

\bibitem{IP98} S. Ivanov, and I. Petrova, {\it Riemannian manifold in which the skew-symmetric curvature operator has pointwise constant
eigenvalues}, Geom. Dedicata {\bf 70} (1998), 269--282.

\bibitem{IZ2005}
S. Ivanov, and S. Zamkovoy, {\it Parahermitian and paraquaternionic
manifolds}, {Differential Geom. Appl.}  {\bf 23} (2005),
205--234.

\bibitem{Kam99} H. Kamada, {\it Neutral hyperk\"ahler structures on primary
Kodaira surfaces}, {Tsukuba J. Math.} {\bf 23} (1999), 321--332.

\bibitem{JK07} J. Kim, {\it On Einstein Hermitian manifolds}, {Monatsh. Math.} {\bf 152} (2007), 251--254.

\bibitem{M03} F. Manhart,
{\it Surfaces with affine rotational symmetry and flat affine metric in $\mathbb{R}^3$},
Studia Sci. Math. Hungar. \textbf{40} (2003), 397--406.

\bibitem{MC05} F. Mart\'in Cabrera, {\it Special almost Hermitian geometry}, {J. Geom. Phys.} {\bf 55} (2005), 450--470.

\bibitem{CS04} F. Mart\'in Cabrera, and A. Swann, {\it Almost Hermitian structures and quaternionic geometries},
{\it Differential Geom. Appl.}   {\bf 21}  (2004), 199--214.

\bibitem{MC06} F. Mart\'in Cabrera, and A. Swann,
{\it Curvature of special almost Hermitian manifolds}, {\it  Pac. J. Math. } {\bf 228} (2006),
165--184.


\bibitem{MS99} A. Mizuhara, and  H. Shima,
{\it Invariant projectively flat connections and its applications},
Lobachevskii J. Math. \textbf{4} (1999), 99--107.

\bibitem{MO08} A. Moroianu, and L. Ornea,
{\it Conformally Einstein products and nearly Kaehler manifolds},
{Ann. Global Anal. Geom.} {\bf 33} (2008), 11--18.
\bibitem{N03}Y.  Nikolayevsky,  {\it Osserman manifolds and Clifford
structures}, {Houston J. Math.} {\bf 29} (2003), 59--75.


\bibitem{N04}
Y. Nikolayevsky, {\it Riemannian manifolds whose curvature operator
$R(X,Y)$ has constant eigenvalues}, {Bull. Austral. Math. Soc.}
 {\bf 70} (2004),  301--319.

\bibitem{Ni1}
Y. Nikolayevsky, {\it Osserman manifolds of dimension $8$}, {
Manuscripta Math.} {\bf 115} (2004), 31--53.



\bibitem{Ni2}
Y. Nikolayevsky, {\it Osserman conjecture in dimension $\neq 8, 16$},
{Math. Ann.} {\bf 331} (2005), 505--522.
\bibitem{O90}
    R. Osserman,   {\it Curvature in the eighties},
    {Amer. Math. Monthly} {\bf  97}, (1990), 731--756.

\bibitem{PSS94} U. Pinkall, A. Schwenk-Schellschmidt, and U. Simon,
{\it Geometric methods for solving Codazzi and Monge-Amp\`ere equations},
Math. Ann. \textbf{298} (1994), 89--100.

\bibitem{PS96}
F. Podesta, and A. Spiro,
{\it Introduzione ai Gruppi di Trasformazioni},
{Volume of the Preprint Series of the Mathematics
Department} V. Volterra of the University of Ancona, Via delle Brecce Bianche, Ancona, ITALY (1996).

\bibitem{Sa04} T. Sato, {\it Examples of Hermitian manifolds with pointwise constant anti-holomorphic sectional
curvature}, {J. Geom.  \bf 80}  (2004), 196--208.

\bibitem{Sa03} T. Sato, {\it Almost Hermitian structures induced from a Kaehler structure which has constant
holomorphic sectional curvature}, {Proc. Amer. Math. Soc.  \bf 131}  (2003),  2903--2909.

\bibitem{SS-62} P. A. Schirokow and A. P. Schirokow,
{\it Affine Differentialgeometrie}, Teubner Leipzig (1962).

\bibitem{R84} {R. Schoen, {\it Conformal deformation of a Riemannian metric to constant scalar curvature},
{J. Differential Geom.} {\bf 20} (1984), {479--495}.}

\bibitem{S95} K. Sekigawa, H. Suga, and L. Vanhecke, {\it Curvature homogeneity for four-dimensional manifolds},  J. Korean Math.
Soc.  {\bf 32}  (1995), 93--101.


\bibitem{S96} K. Sekigawa, {\it Almost Hermitian manifolds satisfying some curvature conditions},  Proceedings of the First International
Workshop on Differential Geometry (Taegu, 1996),  1--11, Kyungpook Natl. Univ., Taegu, 1997.

\bibitem{SSVa92} K. Sekigawa, H. Suga, and L. Vanhecke,
{\it Four-dimensional curvature homogeneous spaces},
{Commentat. Math. Univ. Carol.} {\bf 33} (1992),
261--268.




\bibitem{SSV91} U. Simon, A. Schwenk-Schellschmidt, and H. Viesel,
{\it Introduction to the affine differential geometry of hypersurfaces}, Lecture Notes, Science University of
Tokyo 1991.

\bibitem{Ty74}
W. Takagi,
{\it On curvature homogeneity of Riemannian manifolds},
{T\^ohoku Math. J.} {\bf 26} (1974), 581--585.

\bibitem{S60}
I. M. Singer,
{\it Infinitesimally homogeneous spaces},
{Commun. Pure Appl. Math.} {\bf 13} (1960), 685--697.

\bibitem{ST} I. M. Singer and J. A. Thorpe, {\it The curvature of $4$-dimensional
Einstein spaces}, { 1969
Global Analysis (Papers in Honor of K.
   Kodaira)}, Univ. Tokyo Press, Tokyo, 355--365.

\bibitem{S88} R. Strichartz, {\it Linear algebra of curvature tensors and their covariant
derivatives},
{Can. J. Math}, XL (1988), 1105--1143.

\bibitem{T06} Z. Tang, {\it Curvature and integrability of an almost Hermitian structure},
 {Internat. J. Math.} {\bf 17}, (2006), 97--105.

\bibitem{TV81} F. Tricerri, and L. Vanhecke, {\it Curvature tensors on almost Hermitian manifolds}, {Trans.
Amer. Math. Soc.} {\bf  267}  (1981), 365--397.

\bibitem{TV86} F. Tricerri and L. Vanhecke,
{\it Vari\'et\'es riemanniennes dont le tenseur de courbure est celui d'un espace sym\'etrique
riemannien irr\'eductible},
C. R. Acad. Sci. Paris, S\'er. I {\bf 302} (1986), 233--235.

\bibitem{T68}{N. Trudinger, {\it Remarks concerning the conformal deformation of Riemannian structures on compact manifolds},
{Ann. Scuola Norm. Sup. Pisa} {\bf 22} (1968), 265--274.}

\bibitem{V07} L. Vezzoni, {\it On the Hermitian curvature of symplectic manifolds}, {Adv. Geom.  \bf 7}
(2007), 207--214.


\bibitem{Y60}{H. Yamabe, {\it On a deformation of Riemannian structures on compact manifolds},
{Osaka Math. J.} {\bf 12} (1960), 21--37.}




\end{thebibliography}
\end{document}